\documentclass[11pt]{amsart}
\usepackage{etex}
\pdfoutput=1 
\usepackage{diagbox}
\usepackage{accents}
\usepackage{array}
\usepackage{graphicx, pictex}
\usepackage{stmaryrd}
\usepackage{changepage}
\usepackage{float}
\restylefloat{table}
\usepackage{multirow}
\usepackage[dvipsnames]{xcolor}
\usepackage{amsmath,amssymb,amsthm,mathrsfs}
\allowdisplaybreaks
\def\tildeu#1{\underaccent{\sim}{#1}}

\def\C{{\mathcal{C}}}

\def\d{\Omega}

\def\T{\mathcal{T}}

\newcommand{\tn}[1]{\lVert\kern-1pt\lvert{#1}\rvert\kern-1pt\rVert}
\def\<{{\langle}}
\def\>{{\rangle}}

\def\C{{\mathcal{C}}}

\def\d{\Omega}

\def\d{\Omega}

\def\tb#1{{\|\kern-1pt| #1 \|\kern-1pt|}}
\def\nm2#1#2{\|#1\|_{2,\d_{#2}}}

\def\R{\mathbb{R}}

\def\T{\mathcal{T}}
 \theoremstyle{plain}
 \newtheorem{thm}{Theorem}[section]
 \numberwithin{equation}{section} %% Comment out for sequentially-numbered
 \numberwithin{figure}{section} %% Comment out for sequentially-numbered
 \theoremstyle{plain}
 \newtheorem{prop}[thm]{Proposition}
 \theoremstyle{plain}
 \theoremstyle{plain}
 \newtheorem{theorem}[thm]{Theorem}
 \theoremstyle{plain}
 \newtheorem{corollary}[thm]{Corollary}
\theoremstyle{plain}
 \newtheorem{remark}[thm]{Remark}
 \theoremstyle{plain}
 \newtheorem{lemma}[thm]{Lemma}
\def\C{{\mathcal{C}}}
\def\T{{\mathcal{T}}}

\def\d{{\Omega}}

\def\<{{\langle}}
\def\>{{\rangle}}
\def\R{\mathbb{R}}

\usepackage{amssymb}

\begin{document}

\title[Finite Difference vs Finite Element]
{Connections  Between Finite Difference  and Finite Element Approximations}

\author{Cristina Bacuta}
\address{University of Delaware,
Department of Mathematical Sciences,
501 Ewing Hall 19716}
\email{crbacuta@udel.edu}

\author{Constantin Bacuta}
\address{University of Delaware,
Department of Mathematical Sciences,
501 Ewing Hall, Newark, DE 19716}
\email{bacuta@udel.edu}


\thanks{The work was supported  by NSF-DMS 2011615}%

\begin{abstract}

We present useful connections between  the finite difference and the  finite element methods for a model boundary value problem. We start from the observation that, in the finite element context, the interpolant of the solution in one dimension coincides with the finite element  approximation of the solution.  This result can be viewed as an extension of the Green function  formula for the solution  at the continuous level.  We  write the finite difference and the  finite element  systems such that  the two corresponding linear systems have the same stiffness matrices  and compare the  right hand side load vectors for the two methods.  Using evaluation of the Green function, a formula for the  inverse of the stiffness  matrix is extended to the case of non-uniformly distributed mesh points. We provide an error analysis based on the connection between the two methods, and estimate the energy norm of the difference of the two solutions. Interesting  extensions to the 2D case are  provided.

\end{abstract}
\maketitle

%%%%%%%%%%%%%%%%%            Introduction 
\section{Introduction}

When studying  basic numerical methods for solving boundary value problems (BVP) many resources start with a Finite Difference (FD) approach followed by a separate Finite Element (FE) approach. In this paper, we adopt a new point of view that emphasizes on the connections between the FE and FD methods by solving a standard two point boundary value problem discretized on the same nodes through both  FD and the FE methods. We present the connections that help simplify certain proofs for FD error approximation, and also, that lead to a better understanding of the advantages of each of these two methods.

Consider the two-point boundary value problem 
\begin{equation}\label{2pBVP}
-u''(x) =f(x), \ x \in (0, 1), \ u(0)=u(1)=0.
\end{equation}

It is known that, if $u \in C^2([0, 1])$ is the unique solution of \eqref{2pBVP}, then 
\begin{equation}\label{GF}
u(x) =\int_0^1 G(x,s) f(s)\, ds, \ \text{for all} \ x \in(0,1), 
\end{equation}
where,
\begin{equation}\label{Gfunc}
 G(x,s) =  \left\{
 \begin{array}{rcl}
     s(1-x)  & \mbox{if } \  0\leq s\leq x,\\
      x(1-s)  & \mbox{if } \ x<s\leq 1.\\ 
\end{array} 
\right. 
\end{equation}

 For the discretization of  \eqref{2pBVP}, we divide the interval $[0,1]$ into $n$  subintervals,  using the nodes $0=x_0<x_1<\cdots < x_n=1$ and denote \\ $h_j:=x_j - x_{j-1}, j=1, 2, \cdots, n$. 
 First, consider the finite difference approximation of  $u''(x_j)$  that uses  the quadratic polynomial interpolation  of the solution $u$ at three nodes: $x_{j-1}, x_j$, and $x_{j+1}$:
 \begin{equation}\label{FD1}
 u^{''}(x_j) \approx \frac{2 u(x_{j-1}) }{h_j(h_j +h_{j+1})} + \frac {-2 u(x_j)}{h_j\, h_{j+1}}  + \frac{2 u(x_{j+1})}{h_{j+1}(h_j +h_{j+1})}.
\end{equation}
For the  uniform distribution of the nodes  $x_j=h j$,  $j=0,1, \ldots, n$,  where $h=\frac{1}{n}$,  the approximation \eqref{FD1} becomes the second order  standard centered difference approximation 
\[
 u^{''}(x_j) \approx \frac{u(x_{j-1}) -2 u(x_j) + u(x_{j+1})}{h^2}.
\]

For the general case of non-uniform distributed nodes,  we  let $f_j:=f(x_j)$, and for $u_j\approx u(x_j)$, we solve the system 
  
 \begin{equation}\label{FDs1}
  \left\{
 \begin{array}{rcl}
     u_0 & =0&\\ \\
 \displaystyle  \frac{-2 u_{j-1} }{h_j(h_j +h_{j+1})} + \frac {2 u_j}{h_j\, h_{j+1}}  + \frac{-2 u_{j+1}}{h_{j+1}(h_j +h_{j+1})} & =f_j& \ j=\overline{1,n-1}.\\  \\
     u_n & =0.&
\end{array} 
\right. 
\end{equation}
By multiplying the generic equation in \eqref{FDs1} with $\displaystyle \frac{h_j +h_{j+1}}{2}$,   we get 
 \begin{equation}\label{FDs2}
  \frac{-1}{h_j}\,  u_{j-1} + \left (\frac {1}{h_j} + \frac{1}{h_{j+1}} \right )  u_j + \frac{-1} {h_{j+1}} u_{j+1} = \frac{h_j +h_{j+1}}{2}\,  f_j.
\end{equation}
 
% \begin{equation}\label{Gfunc}
%  \left\{
% \begin{array}{rcl}
%     u_0 & =0&\\
 %    -\frac{u_{j-1} -2 u_j +u_{j+1}}{h^2} & =f_j& \ j=1,2,\cdots,n-1.\\ 
 %    u_n & =0.&
%\end{array} 
%\right. 
%\end{equation}
Denote ${u}^{FD} := [ u_1, u_2, \cdots, u_{n-1}]^T$, and  $\tilde{f}:=[ f_1, f_2, \cdots, f_{n-1}]^T$. Let  $W$  be the $(n-1)\times(n-1)$ diagonal matrix with entries: $\{\frac{h_1 +h_{2}}{2}, \cdots, \frac{h_{n-1}+h_{n}}{2} \}$ and let  $S$  be the stiffness
tridiagonal $(n-1)\times(n-1)$ matrix
\begin{equation} \label{S} 
S=\left [\begin{matrix}
 \frac{1}{h_1}+ \frac{1}{h_2},& -\frac{1}{h_2} &   & & & \\
-\frac{1}{h_2} &  \frac{1}{h_2}+ \frac{1}{h_3} & -\frac{1}{h_3}  &  & &  \\
% & -\frac{1}{h_3} &  \frac{1}{h_3}+ \frac{1}{h_4} & -\frac{1}{h_4}  &  &   \\
  \   & \ddots & \ddots  & \ddots  & &  \\
    & & & & & \\
  &  & & -\frac{1}{h_{n-2}}  &  \frac{1}{h_{n-2}}+ \frac{1}{h_{n-1}} &   -\frac{1}{h_{n}} \\
 & &  & & -\frac{1}{h_{n-1}}  &  \frac{1}{h_{n-1}}+ \frac{1}{h_n}  \\
\end{matrix} \right].
\end{equation}

Then, using \eqref{FDs2}  and the fact that $u_0=u_n=0$, the system \eqref{FDs1} is equivalent to  
 \begin{equation}\label{FDs3}
 S\, {u}^{FD} = W\, \tilde{f}.
\end{equation}

%For the uniform distribution of the nodes   $x_0, x_1, \cdots, x_n$, 
%denoting by $A:= $ the  following $(n-1)\times (n-1)$ {\it  tridiagonal} matrix  $A=tridiag(-1, 2,-1)$ or 
%$$ A=\left [\begin{matrix}
%\ 2 & -1 &  & & & \\
%-1 & \ 2 & -1 &  & & \\
%\  &  -1 & \ 2 & -1 &  &  \\
% \  & & \ddots & \ddots  & \ddots  &   \\
% &  & & -1 & \ 2 & -1 \\
% & &  & & -1 & \ 2  \\
%\end{matrix} \right],$$
%\vspace{0.2in}
%we have that $S= \frac{1}{h} \, A$ and $W= h I$, with $I$ the identity matrix on $\R^{n-1}$. Thus, the system to solve for $ {u}^{FD}$    is the  \begin{equation}\label{FDs4}
 %\frac{1}{h} \, A \, {u}^{FD} = h\, \tilde{f}.
%\end{equation}

For a better comparison of the finite difference  solution $ {u}^{FD}$ to the finite element  approximation $ {u}^{FE}$, we prefer to write the linear system in the form \eqref{FDs3}.

Secondly, for the FE discretization, we  will use the standard notation
\[ 
a(u, v) = \int_0^1 u'(x) v'(x) \, dx, \ \text{and} \ (f, v) = \int_0^1  f(x) v(x) \, dx.
 \]
Then, the variational formulation of \eqref{2pBVP} is : Find $u \in V:= H^1_0(0,1)$ such that
 \begin{equation}\label{VF}
a(u, v) = (f, v), \ \text{for all} \ v \in V.
\end{equation}

 Also, we consider the mesh nodes  on $[0,1]$ as being the same as those used for the finite difference discretization:   $0=x_0<x_1<\cdots < x_n=1$,  and define  the corresponding   discrete space  $V_h$  as  the subspace of $ V = H^1_0(0,1)$, given by
 \[ 
 V_h = \{ v_h \in V \mid v_h \text{ is linear on each } [x_j, x_{j + 1}]\},
 \]
  i.e., $V_h$ is the space of all {\it piecewise linear continuous functions} with respect to the given nodes that are zero at $x=0$ and $x=1$.  We consider the nodal basis $\{ \varphi_j\}_{j = 1}^{n-1} \subset V_h$ such that $\varphi_i(x_j ) = \delta_{ij}$. 
 Thus, the discrete  variational formulation of \eqref{2pBVP} is: Find $u_h \in V:= H^1_0(0,1)$ such that
 \begin{equation}\label{dVF}
a(u_h, v_h) = (f, v_h), \ \text{for all} \ v_h \in V_h.
\end{equation}
We look for  $u_h=u_h^{FE}$ with  the nodal basis expansion
 \[
 u_h := \sum_{i=1}^{n-1} u_i \varphi_i, \ \text{where} \ \ u_i=u_h(x_i).
 \]
If we consider the test functions $v_h=\varphi_j, j=1,2,\cdots,n-1$ in \eqref{dVF}, we  obtain   the system 
\begin{equation}\label{eq:FEsol}
S\, u^{FE} = \tildeu{f},
\end{equation}
for finding the coefficient vector  $u^{FE}=[u_1, \cdots, u_{n-1}]^T$, where \\ $\tildeu{f}:=[(f, \varphi_1), \cdots (f, \varphi_{n-1})]^T$ and  $S$ is the $(n-1)\times (n-1)$ {\it tridiagonal} matrix with entries $S_{ij}=a(\varphi_i', \varphi_j')$. 
It is an easy exercise to verify that $S$ defined for the finite element discretization, coincides with the matrix defined in \eqref{S} for the finite difference discretization. 

We note that  for any $\alpha \in \R^{n-1}, \alpha \neq 0$ 
\[
\alpha^T S \, \alpha = a(v_h, v_h) =  \int_0^1 (v'_h(x))^2 \, dx>0,
\]
where $v_h = \sum_{i=1}^{n-1} \alpha_i \varphi_i$. Consequently, the stifnees matrix $S$ is invertible.

We can also define the $(n-1)\times(n-1)$ Green matrix 
\[
\tilde{G} =\left [G(x_i,x_j) \right ]_{i,j=\overline{1,n-1}}.
\]

The purpose of the paper is to show how the  two types of discretization are connected. 
%and provide an estimate for the  error $u^{FD} -u^{FE}$.   
Using  that  the Green function \eqref{Gfunc} can be viewed  as a scaled  finite element  function, we will prove that  the common stiffness matrix $S$ satisfies  $ S^{-1}=\tilde{G}$ in the general non-uniform case. While the identity is well known for the uniform distributed nodes, see Section 12.2.2 in \cite{QSS}, by the best knowledge of the authors, the identity seems to  be not available  for the non-uniform case or in the finite element context. We  will also show that  the right hand side vector of  \eqref{FDs3}, that defines the finite difference solution,  is obtained by the Composite Trapezoid Rule (CTR) approximation of  the entries  of the dual vector  $\tildeu{f}$.  As a consequence, we can  compare the $u^{FD}$  and $u^{FE}$  and  estimate  the energy norm of the  error $u^{FD} -u^{FE}$.

The rest of the paper is organized as follows: In Section \ref{sec:FDvsFE}, we  present the similarity between the  finite element and the finite difference  linear systems   for  the model problem \eqref{2pBVP}. In Section \ref{sec:error}, we use the formulas for the FD and the FE solutions in order to compare the two discretizations and  do an error analysis.  Two interesting extensions to the 2D case are presented in Section \ref{sec:2D}. %We conclude with Section \ref{sec:conclusion}. 

%%%%%%%%%%%%%%%%%%%%%%%%%%%%%%%%%%%
\section{The connections between the finite element and the finite difference of the 1D model problem}\label{sec:FDvsFE}
%%%%%%%%%%%%%%%%%%%%%%%%%%%%%%%%%%%
In this section we investigate the systems for the FE and FD discretizations of \eqref{2pBVP} by relating the matrices and the load vectors for the two systems. The following lemma gives a formula for the general  component of the dual vector $\tildeu{f}$ and the evaluation of the solution $u$  at the neighboring nodes viewed as degrees of freedom for the finite element approximation. 
In addition, the formula connects with the FD discretization via the formula  \eqref{FDs2}. 
%%%%%%%%%%%%%%%%%%%%%%%%%%%%%%%%%%%%
\subsection{A formula for the finite element  solution}
%%%%%%%%%%%%%%%%%%%%%%%%%%%%%%%%%%%%
\begin{lemma}\label{fvarphi}
Let $u$ be the solution of \eqref{2pBVP}, and $\{\varphi_j\}_{j=1,2,\cdots, n-1}$ be the nodal basis for $V_h$. Then, for $j=1,2,\cdots,n-1$, we have

\begin{equation}\label{lemma1fphij}
(f,\varphi_j)=\displaystyle -\frac{1}{h_j} u(x_{j-1})+\left(\frac{1}{h_j}+\frac{1}{h_{j+1}}\right) u(x_j)-\frac{1}{h_{j+1}} u(x_{j+1})
\end{equation}

\end{lemma}

\begin{proof}
From the variational formulation \eqref{dVF}, we have $(f,\varphi_j)=a(u,\varphi_j)=\displaystyle \int_0^1u'(x)\varphi'_j(x) dx=\frac{1}{h_j}\int_{x_{j-1}}^{x_j} u'(x) dx+\frac{-1}{h_{j+1}}\int_{x_{j}}^{x_{j+1}} u'(x) dx=-\frac{1}{h_j} u(x_{j-1})+\left(\frac{1}{h_j}+\frac{1}{h_{j+1}}\right) u(x_j)-\frac{1}{h_{j+1}} u(x_{j+1})$.
\end{proof}

Next, for any $s\in (0,1)$, we define the function
\begin{equation}
 \phi_s(x) =  \left\{
 \begin{array}{ccl}
    \displaystyle  \frac{x}{s}  & \mbox{if } \  0\leq x\leq s,\\ \\
    \displaystyle   \frac{1-x}{1-s}  & \mbox{if } \ s<x\leq 1.\\ 
\end{array} 
\right. 
\end{equation}

Note that, if $s=x_j$, $j=1,2,\cdots,n-1$, then $\phi_{x_j}$ is a piecewise linear function that belongs to $V_h$.

\begin{corollary} (Green's Formula)
For any $s\in (0,1)$, we have

\begin{equation}\label{fBigphij}
(f,\phi_s)=\displaystyle \frac{1}{s(1-s)} u(s)
\end{equation}

\end{corollary}
\begin{proof} In the Lemma \eqref{fvarphi}, take $n=2$, $x_0=0$, $x_2=1$, $x_1=s$, $h_1=s$, $h_2=1-s$.
\end{proof}

Note that the equation  \eqref{fBigphij} is in fact the Green's formula \eqref{dVF}. This is because 

\begin{equation}\label{GBigphis}
G(x,s)=s(1-s) \phi_s(x)=G(s,x), \ \text{for all}  \ x, s \in (0, 1).
\end{equation}

This formula shows that the Green's Function can be interpreted as a scaled finite element function.

\begin{corollary}\label{uFF} Let $\displaystyle u_h := \sum_{i=1}^{n-1} u_i \varphi_i$  be the finite element solution of  \eqref{dVF}, and let ${u}^{FE} := [ u_1, u_2, \cdots, u_{n-1}]^T$. Then $u_h$ coincides with the linear interpolant of $u$ on the nodes $x_0,x_1, \ldots, x_n$. 

In other words, $u_j=u(x_j)$, $j=1,2,\cdots,n-1$.

\end{corollary}
\begin{proof} Since the stiffness matrix $S$ for the system \eqref{dVF} is given by \eqref{S}, from the equation \eqref{lemma1fphij}  we have that $\tilde{u}=[u(x_1), u(x_2), \ldots, u(x_{n-1})]^T$ solves $S\tilde{u}=\tildeu{f}$. On the other hand, from \eqref{eq:FEsol} we have that 
$S\, u^{FE}=\tildeu{f}$. Using  that the matrix $S$ is invertible, we conclude that $u^{FE}=\tilde{u}$.
\end{proof}

\begin{corollary} If $S$ is the stiffness matrix defined in \eqref{S}, then 
\begin{equation}\label{SinverseG}
S^{-1}=\tilde{G},
\end{equation}
consequently, $u^{FE}=S^{-1}\tildeu{f}= \tilde{G} \tildeu{f}$. 
\end{corollary}

\begin{proof}
Apply \eqref{fBigphij} for $s=x_j$, $j=1,2,\cdots,n-1$ to get 
\begin{equation}\label{uxPhij}
u(x_j)=x_j(1-x_j)\int_0^1\phi_{x_j}(x)f(x) dx.
\end{equation}
Obviously, by the definition of $\phi_{x_j}$ we have $\phi_{x_j}\in V_h$. Hence, $\phi_{x_j}$ coincides with its interpolant on the nodes $x_0,x_1, \ldots, x_n$. Thus,
\begin{equation}\label{Phivarphi}
\phi_x(x)=\sum_{i=1}^{n-1} \phi_{x_j}(x_i) \varphi_i(x).
\end{equation}

For the equations \eqref{uxPhij} and \eqref{Phivarphi}, we get

\begin{equation}\label{eq:FESf}
u(x_j)=x_j(1-x_j)\sum_{i=1}^{n-1} \phi_{x_j}(x_i) \int_0^1f(x)\varphi_i(x) dx.
\end{equation}

 Since $x_j(1-x_j)\phi_{x_j}(x_i)=G(x_j,x_i)$, we obtain $\displaystyle u(x_j)=\sum_{i=1}^{n-1} G(x_j,x_i)(f,\varphi_i)$, $j=1,2,\cdots,n-1$, which combined with the fact that the finite element solution coincides with the interpolant as given by Corollary \ref{uFF}, leads to $u^{FE}=\tilde{G}\tildeu{f}$. Since $u^{FE}=S^{-1}\tildeu{f}$ for any $\tildeu{f} \in \R^{n-1}$, we conclude that the formula \eqref{SinverseG} holds.
\end{proof}
In order to obtain the formula \eqref{eq:FESf}, it is essential that $\phi_{x_j}\in V_h$, and that  $\phi_{x_j}$ can be expanded in the basis $\{ \varphi_i\}$. 

%%%%%%%%%%%%%%%%%%%%%%%%%%%%%%%%%%%%
\subsection{A formula for the finite difference solution}
%%%%%%%%%%%%%%%%%%%%%%%%%%%%%%%%%%%%

Let us consider the finite difference discretization \eqref{FDs1} that is equivalent to \eqref{FDs2} and  builds the system \eqref{FDs3}.

For a continuous function $\theta : [0,1]\to \R$ such that $\theta(0)=\theta(1)=0$, the composite trapezoid rule (CTR) on the nodes $x_0,x_1, \ldots, x_n$ is 
\begin{equation}\label{CTR}
\int_0^1\theta(x)dx \approx T_n(\theta):=\sum_{i=1}^{n-1} \theta(x_i) \frac{h_i+h_{i+1}}{2}.
\end{equation}
Note that if $\displaystyle I(\theta)=\sum_{i=1}^{n-1} \theta(x_i)\, \varphi_i$ is the interpolant of $\theta$, and by using  that $\int_0^1 \varphi_i(x)dx=\displaystyle \frac{h_i+h_{i+1}}{2}$, then the CTR formula  becomes
\[
T_n(\theta)=\int_0^1I(\theta)(x)dx. 
\]
We observe that  the diagonal entries of  $W$ in equation \eqref{FDs3} are  exactly  the weights of the quadrature formula \eqref{CTR}. 

Next, we will show that, componentwise, the finite difference solution $u_j$  is the CTR approximation of the solution $u(x_i)$ given by the Green's formula:
\begin{equation}\label{Guj}
u(x_j)=\int_0^1 G(x_j,x)f(x)dx.
\end{equation} 

For each $j=1,2,\cdots,n-1$, we define 
\[
\theta_j(x)=G(x_j,x)f(x), \ x \in [0,1], \ \text{and}
\]
\[
\begin{aligned}
w_j & :=T_n(\theta_j)=\displaystyle \sum_{i=1}^{n-1}\frac{h_i+h_{i+1}}{2}\theta_j(x_i)
=\displaystyle \sum_{i=1}^{n-1} G(x_j,x_i)\frac{h_i+h_{i+1}}{2}f(x_i) \\
& =(\tilde{G}W\tilde{f})_j.
\end{aligned}
\]
Consequently, we have
\begin{equation}\label{eq:w}
w= \tilde{G}W\tilde{f}.
\end{equation}
\begin{prop}\label{WT}
The finite difference solution $u^{FD}$ of \eqref{FDs3} coincides with the vector $w=[w_1, w_2, \ldots, w_{n-1}]^T$.
\end{prop}
\begin{proof}
From  the formula  \eqref{SinverseG}, we have  that $u^{FD}=\tilde{G} W\tilde{f}$ which compared to  \eqref{eq:w}, leads to $u^{FD}=w$.
\end{proof}

\begin{remark} We can prove the Proposition \ref{WT} without using finite element arguments needed to prove  \eqref{SinverseG}. We can just verify that the vector $w$ satisfies the system $Sw=W\tilde{f}$ by essentially using  that $\phi_{x_j}$  is linear on the intervals $[0, x_j]$ and $[x_j, 1]$.
\end{remark}
%%%%%%%%%%%%%%
\subsection{Comparison}
%%%%%%%%%%%%%%%%%%
 By rewriting the FD system, we have that the  FD and the FE linear systems have the same stiffness matrices.  The corresponding solutions are given by  
\[
u^{FE}=S^{-1}\tildeu{f}, \ \text{and} \  u^{FD}=S^{-1}W\tilde{f},\ \text{where} 
\]
\[
(\tildeu{f})_j=(f,\varphi_j)=\int_0^1 f(x)\varphi_j(x)\, dx,  \ \text{and} \ (W\tilde{f})_j=\frac{h_j+h_{j+1}}{2}f(x_j) . 
\]
The right hand sides of  the two systems are componentwise related by 
\begin{equation}\label{fTn}
(f,\varphi_j)=T_n(f\varphi_j).
\end{equation}

%%%%%%%%%%%%%%%%%%%%%%%%%%%%%%%%%%%%%
\section{Error Estimates} \label{sec:error}
%%%%%%%%%%%%%%%%%%%
 As a consequence of the previous section, we provide error estimates for $u -u_h^{FD}$  in both the discrete infinity norm and the energy norm, and estimate the energy norm of the difference $u_h^{FE} -u_h^{FD}$. 

%%%%%%%%%%%%%%%%%%%%%%%%%%%%%%%%%%%%%
\subsection{The infinity norm estimate for the FD solution} 
%%%%%%%%%%%%%%%%%%%
From Lemma \ref{WT}, we can estimate the error $$
\displaystyle  \max_{j=\overline{1,n-1}} |u(x_j)-u_j^{FD}| $$  in the general non-uniform case.

For the uniform case, it is well known (see e.g.  \cite{QSS}), that  if $ f\in \C^2([0,1])$, then
\begin{equation}\label{maxerrc}
\displaystyle \max_{1\leq j\leq n} |u(x_j)-u_j^{FD}|\leq \frac{h^2}{96}\|f''\|_{\infty},
\end{equation}
where for a continuous function $\theta$ on $[a,b]$, we define $\displaystyle \|\theta\|_{\infty, [a, b]}:=\max_{x\in[a,b]}|\theta(x)|$ and 
$\|\theta\|_{\infty}:=\|\theta\|_{\infty, [0, 1]}$. 

For  extending \eqref{maxerrc} to the general nonuniform case,  we will use the known  formula for the CTR error as follows: 
For a function $\theta:[a,b]\to \R$, \\ $\theta \in \C^2([a,b])$, and the nodes $a\leq x_0 <x-1<\ldots <x_n=b$, we have
\begin{equation}\label{CTRE}
\displaystyle \left|\int_a^b\theta(x)dx - T_{n,[a,b]}\right|\leq\frac{b-a}{12} h^2\|\theta''\|_{\infty,[a,b]}, 
\end{equation}
where $\displaystyle T_{n,[a,b]}(\theta)=\sum_{i=1}^n (x_i-x_{i-1})\frac{\theta(x_i)+\theta(x_{i-1})}{2}$ and $h=\displaystyle \max_{i=\overline{1,n}}(x_i-x_{i-1})$.

\begin{theorem} Let $f\in  \C^2[0,1])$, and  let  $u$ be the solution of the boundary value problem \eqref{2pBVP}, and $u^{FD}=[u_1,u_2,\ldots, u_{n-1}]^T$ be the finite difference solution of \eqref{FDs3}. If $h_j=x_j-x_{j-1}$ for $j=1,2,\ldots, n$ and $\displaystyle h:=\max_{j=\overline{1,n}} |h_j|$, then, we have
\begin{equation}\label{EDEinf}
\displaystyle \max_{1\leq j \leq n-1}|u(x_j)-u_j|\leq \frac{h^2}{48}(4\|f'\|_{\infty}+\|f''\|_{\infty}).
\end{equation}
\end{theorem}
\begin{proof}
From the proof of Corollary \ref{SinverseG} and the equation \eqref{uxPhij}, we have 
\[
u(x_j)=x_j(1-x_j)\int_0^1 \phi_{x_j}(x) f(x)dx. 
\]
From  Proposition  \ref{WT}, equations \eqref{CTR} and \eqref{fTn}, we also have that 

\begin{equation}\label{uFDf}
\begin{aligned}
u_j&=(\tilde{G}W\tilde{f})_j= \displaystyle \sum_{i=1}^{n-1} G(x_j,x_i)\frac{h_i+h_{i+1}}{2}f(x_i) \\ 
&= \displaystyle x_j(1-x_j)  \sum_{i=1}^{n-1} \phi_{x_j}(x_i)\frac{h_i+h_{i+1}}{2} = \displaystyle x_j(1-x_j) T_n(\phi_{x_j}f).
\end{aligned}
\end{equation}
Consequently, the error  $E_j=u(x_j)-u_j$  satisfies 
\begin{equation}
\begin{aligned}
E_j&= \displaystyle x_j(1-x_j)\left( \int_0^1\phi_{x_j}(x)f(x)dx- T_n(\phi_{x_j}f)\right)=\\ 
&=x_j(1-x_j)\left(\int_0^{x_j} \phi_{x_j}(x)f(x)dx - T_{n, [0, x_j]}(\phi_{x_j}f)\right) +\\
& + x_j(1-x_j)\left(\int_{x_j}^1 \phi_{x_j}(x)f(x)dx- T_{n, [ x_j,1]}(\phi_{x_j}f)\right).
\end{aligned}
\end{equation}

By noting that $\phi_xf$ is $\C^2$ on $[0,x_j]$ and $[x_j,1]$, applying the estimate \eqref{CTRE} on the intervals $[0,x_j]$ and $[x_j,1]$  for $\theta(x)=\phi_{x_j}(x)f(x)$, and using
\[
\displaystyle \max_{x\in [0,x_j]}|\theta''(x)|\leq \frac{2}{x_j}\|f'\|_{\infty}+\|f''\|_{\infty}, \ \text{and}
\]
\[
\displaystyle \max_{x \in [x_j,1]}|\theta''(x)|\leq \frac{2}{1-x_j}\|f'\|_{\infty}+\|f''\|_{\infty},
\]
we obtain 
\begin{equation*}
\begin{aligned}
\displaystyle |E_j| & \leq  x_j(1-x_j) \frac{h^2}{12}x_j\left(\frac{2}{x_j}\|f'\|_{\infty,[0,x_j]}  
 +\|f''\|_{\infty}\right)+ \\
&+x_j(1-x_j) \frac{h^2}{12}(1-x_j)\left(\frac{2}{1-x_j}\|f'\|_{\infty,[x_j,1]}+\|f''\|_{\infty}\right).
\end{aligned} 
\end{equation*}
Since  $x_j(1-x_j) \leq \frac{1}{4}$, the above estimate  leads to \eqref{EDEinf}.
\end{proof}
\begin{remark} It is known that \eqref{FD1} is only first order accurate if $h_j\neq h_{j+1}$,  see  e.g., Section 1.3 in \cite{LeVeque07}. The estimate \eqref{EDEinf} shows that  the FD method based on \eqref{FDs1} or \eqref{FDs2}  is globally $O(h^2)$.
We also emphasize that our  proof is done using the formula of the solution \eqref{uFDf}, and it is not based on \\  $\delta$-functional arguments as usually done for proving \eqref{maxerrc}, see e.g., \cite{QSS}. 
\end{remark}

 Based on the CTR approximation, the formula \eqref{uFDf} allows for further error analysis in  the energy norm. 

%%%%%%%%%%%%%%%%%%
\subsection{Energy Norm  Errors}
%%%%%%%%%%%%%%%%%%%%%

For $v\in V=H_0^1(0, 1)$, define the energy norm
 \begin{equation}
|v|=|v|_a=\left(\int_0^1 (v'(x))^2dx\right)^{\frac{1}{2}}=(a(v,v))^{\frac{1}{2}}.
\end{equation}

If the vector $u^{FD}=[u_1, \ldots, u_{n-1}]^T$ is the solution of the finite difference system \eqref{FDs3}, then we define the corresponding function in $V_h$ as 

$u_h^{FD}:=\displaystyle \sum_{i=1}^{n-1} u_i\varphi_i$. From the previous sections, we know that $u^{FD}$ satisfies the equation $Su^{FD}=\tilde{G}(W\tilde{f})$, where $W\tilde{f}$ is defined by \eqref{S} and \eqref{FDs3}, and $\displaystyle \frac{h_i+h_{i+1}}{2}f(x_i)=T_n(f\varphi_j)$. Thus, since the finite element stiffness matrix is still $S$, we have that 
\[
a(u_h^{FD}, \varphi_j)=T_n(f\varphi_j), \ \ \text{for all} \ j=1,2,\ldots, n-1.
\]
Using the linearity of $a(u_h^{FD},\cdot)$ and $T_n(f,\cdot)$, we have that 
\begin{equation}\label{uhFD}
a(u_h^{FD},v_h)=T_n(fv_h),  \ \text{for all} \ v_h \in V_h.
\end{equation}

On the other hand, we have 
\begin{equation}\label{uhFE}
a(u_h^{FE}, v_h)=(f,v_h)=\int_0^1fv_h \ dx,  \ \text{for all} \ v_h\in V_h.
\end{equation}

From the equations \eqref{uhFD} and \eqref{uhFE}, we obtain 
\[
a(u_h^{FE}-u_h^{FD}, v_h)=F_h(v_h):=\int_0^1fv_h \ dx-T_n(fv_h). 
\]
Consequently, 
\begin{equation}\label{Fh}
|u_h^{FE}-u_h^{FD}|=\sup_{v_h \in V_h} \frac{F_h(v_h)}{|v_h|}:=\|F_h\|_{V_h^*}
\end{equation}
where the supremum is taken over all non-zero vectors. 
\begin{theorem}\label{th:FED}
Assume that $f\in\C^2([0,1])$, and that $u_h^{FE}$ and $u_h^{FD}$ are the finite element and the finite difference corresponding solutions of the boundary value problem \eqref{2pBVP} on the nodes $x_0=0<x_1<\ldots<x_n=1$.

 Let $\displaystyle h=\max_{i=\overline{1,n}} (x_i-x_{i-1})=\displaystyle \max_{i=\overline{1,n}}(h_i)$. Then, 
\begin{equation}\label{udiff}
|u_h^{FE}-u_h^{FD}|\leq h^2\left( \frac{\|f''\|_{\infty}}{12}+\frac{\|f'\|_{\infty}}{6}\right).
\end{equation}
\end{theorem}
\begin{proof}
According to the equation \eqref{Fh}, we just have to find an upper bound for $\|F_h\|_{V_h^{\star}}$. Let $\theta_{i}(x):=f(x)v_h(x)$, for $x\in [x_{i-1},x_i]$, $i=1,2,\ldots, n$. Using the trapezoid rule error formula for $\theta_i$, we  have that 
\begin{equation}\label{Fvh}
F_h(v_h)=-\frac{1}{12}\sum_{i=1}^n h_i^3 \theta_i ''(\xi_i),  \ \xi_i \in [x_{i-1},x_i].
\end{equation}
To simplify the notation, we assume that $\displaystyle v_h=\sum_{i=1}^{n-1} \alpha_i \varphi_i$. Then, on the interval $(x_{i-1},x_i)$, $v'_h= \frac{\alpha_i - \alpha_{i-1}}{h_i}$, and  from the equation \eqref{Fvh}, we have 
\[
F_h(v_h)=\displaystyle -\frac{1}{12}\left(\sum_{i=1}^n  h_i^3 f''(\xi_i)v_h(\xi_i)+2\sum_{i=1}^n h_i^3\,  \frac{\alpha_i-\alpha_{i-1}}{h_i}f'(\xi_i) \right). 
\]
We have  $h_i\leq h$ and $|v_h(x)| \leq \|v_h\|_{\infty}\leq |v_h|$, for all $v_h \in V_h$. Thus, 
\begin{equation}\label{Fh2}
|F_h(v_h)| \leq \frac{1}{12}h^2|v_h|\sum_{i=1}^n h_i |f''(\xi_i)|+\frac{1}{6} h^2 \|f'\|_{\infty}\sum_{i=1}^n|\alpha_i - \alpha_{i-1}|. 
\end{equation}
Using the discrete  mean value theorem, (see e.g., Theorem 9.1  in  \cite{QSS}), for some $\xi \in (0, 1)$, we have
\[
\sum_{i=1}^n h_i |f''(\xi_i)|=|f''(\xi)| \sum_{i=1}^n h_i\leq \|f''\|_{\infty}.
\]
For the second sum in \eqref{Fh2}, the  Cauchy-Schwarz inequality gives
\[
\begin{aligned}
 \sum_{i=1}^n|\alpha_i-\alpha_{i-1}|&  =\sum_{i=1}^n h_i^{\frac{1}{2}} \frac{|\alpha_i-\alpha_{i-1}|}{h_i^{\frac{1}{2}}} 
  \leq \left(\sum_{i=1}^n h_i\right)^{1/2}\left(\sum_{i=1}^n  \frac{|\alpha_i-\alpha_{i-1}|^2}{h_i}\right)^{1/2}\\ &=\left(\sum_{i=1}^n \int_{x_{i-1}}^{x_i} (v'_h)^2 dx\right )^{1/2}=|v_h|.
 \end{aligned}
 \]
Combining the last two estimates with  \eqref{Fh2}, gives 
\[
|F_h(v_h)|\leq \frac{1}{12} h^2 |v_h|  \|f''\|_{\infty}+\frac{1}{6}h^2\|f'\|_{\infty} |v_h|.
\] 
Hence, 
\[
\|F_h\|_{V_h^*} \leq h^2\left( \frac{\|f''\|_{\infty}}{12}+\frac{\|f'\|_{\infty}}{6}\right), 
\]
which, together with \eqref{Fh}  proves the statement of the theorem.
\end{proof}
\begin{corollary}
Under the same assumptions as in the Theorem \ref{th:FED}, we have  
\begin{equation}\label{ChFD}
|u-u_h^{FD} |\leq Ch(\|f''\|_\infty + \|f\|_{L^2(0, 1)}),
\end{equation}
 where $C$ is a constant independent of $h$ or $f$.
\end{corollary}
\begin{proof} It is known that the following estimate for the finite element error in the energy norm holds

\begin{equation}\label{ChFE}
|u-u_h^{FE} |\leq Ch \|f\|_{L^2(0, 1)}.
\end{equation}
 where $C$ is a constant independent of $h$ or $f$, see e.g. \cite{braess, brenner-scott,  ern-guermond}.

Now, the inequality \eqref{ChFD} follows  using \eqref{ChFE}, \eqref{udiff}, and  the triangle inequality
\[
|u-u_h^{FD}| \leq |u-u_h^{FE}|  + |u_h^{FE}-u_h^{FD}|.
\]

\end{proof}

%%%%%%%%%%%%%%
\section{Interesting extensions to the 2D case}\label{sec:2D}
%%%%%%%%%%%%%%

Let $\Omega\subset\R^2$ be a polygonal domain, and $f \in L^2(\Omega)$ or $f$ is a continuous functional on $H^1_0(\Omega)$. Consider the model problem: Find $u$  such that 

\begin{equation}\label{2DBVP}
 \left\{
 \begin{array}{rccl}
     -\Delta u = & f  & \text{in } \ \Omega,\\
      u =& 0 & \text{on} \  \partial \Omega.
\end{array}  
\right. 
\end{equation}

The corresponding variational or weak formulation of \eqref{2DBVP} is: 

Find $u \in H^1_0(\Omega)$ such that 
\begin{equation}\label{2Dvf}
a(u,v) = (f,v) ,\ \text{for all} \ v \in  H^1_0(\Omega),
\end{equation}
where, for any $u,v\in  H^1_0(\Omega)$, and for  $f \in L^2(\Omega)$, 
\[
a(u,v): =\int_\Omega \nabla u\cdot \nabla v \, dx, \ \text{and}\  (f,v):=\int_\Omega f(x) v(x) \, dx.
\]
The existence and uniqueness of the solution of \eqref{2Dvf} is well know, see e.g., \cite{braess, brenner-scott,  ern-guermond, demko-oden}.
For the discretization of \eqref{2DBVP}, we let $\T_h$ be a triangulation of $\Omega$ and consider $V_h$ a conforming finite element space with a nodal basis $\{\varphi_1, \cdots, \varphi_n\}$ associated with the mesh $\T_h$. 
If $\varphi_j$ is a generic basis function with support $D_j \subset \Omega$, we assume
\begin{equation*}
\displaystyle D_j= \bigcup_{i=1}^{n_j} T_i, \ \text{with}  \ T_i\ \text{being a triangle in} \ \T_h. 
\end{equation*}
The formula \eqref{lemma1fphij} is a formula of significant importance because it shows the strong connection between the FE and the FD discretizations  in 1D. 
A corresponding formula for the 2D case can be obtained by applying the Green's formula for the function $f\varphi_j = -\Delta u \varphi_j$ to get

\begin{equation}\label{2DF}
\begin{aligned}
(f,\varphi_j) &  = \int_\Omega  -\Delta u \varphi_j \, dx = \int_\Omega \nabla u\cdot \nabla \varphi_j\, dx =\\
                  & = -\sum_{T_i\subset D_j}\int_{T_i} (\Delta \varphi_j)\, u\, dx +\sum_{T_i\subset D_j} \int_{\partial T_i} (\nabla u\cdot {\bf n}) \, u \, ds.
\end{aligned}
\end{equation}
where ${\bf{n}}$ is the outer normal vector to the integration domain's boundary. 
%\left (\frac{\partial \varphi_j}{\partial n}\right )
Similary to the 1D case,  the components $(f,\varphi_j)$ of the dual vector $f$ relative to   the basis $\{\varphi_1, \cdots, \varphi_n\}$  can be written as linear combinations of degrees of freedom acting on $u$. However, in this case, the degrees of freedom have to involve integrations on triangles and on the edges of $\T_h$.  We can redefine an  interpolant using the new degrees of freedom in order to have the FE solution agree with the interpolant.
However,  the main challenge  in this case is to match the number of degrees of freedom with the number $n$ of basis functions when considering the system given by \eqref{2DF}. Except for very simple meshes, this seems to be difficult to achieve. 

Next, we consider  two special 2D cases, leaving the challenge to extend this idea to the  more general  cases of meshes in 2D or 3D for  future work.

%%%%%%%%%%%%%%
\subsection{The Green function  for the 1D BVP as a 2D discretization function fo the 2D BVP}

Consider the problem \eqref{2DBVP} on the unit square $\Omega=(0, 1) \times (0, 1)$, and define the mesh $\T_1=T_1\cup T_2$, 
where $T_1, T_2$ are the two triangles determined by the positive slope diagonal $\Gamma$  of $\Omega$. We define $V_h=V_1= span\{G\}$, where $G=G(x,y)$ is the function defined in \eqref{Gfunc}. 

We have  $G\in H^1_0(\Omega)$, and if $u$ is the solution of 
\eqref{2DBVP}, then 
\begin{equation}\label{2DG}
\begin{aligned}
(f,G) &  = \int_\Omega  -\Delta u\, G\, dxdy = \int_\Omega \nabla u\cdot \nabla G\, dxdy =\\
                  & = \int_{T_1} -\Delta G \, u\, dxdy  \int_{\partial T_1} (\nabla G\cdot {\bf n}_{T_1})\, u \, ds +\\
                  & + \int_{T_2} -\Delta G \, u\, dxdy  \int_{\partial T_2} (\nabla G\cdot {\bf n}_{T_2})\, u \, ds.
\end{aligned}
\end{equation}
Using that $\Delta G= 0 \ \ \text{on} \ \Omega \backslash \Gamma$ and on $\Gamma= \{(x,y) \in\Omega |\  x=y\}$, and that
\[
\nabla G _{|_{T_1}} \cdot {\bf n}_{T_1}+ \nabla G _{|_{T_2}} \cdot {\bf n}_{T_2}=\sqrt{2},
\]
from \eqref{2DG}, we obtain a  formula for $\int_\Gamma u\, ds $. 
%\frac{\partial G_{|_{T_1}}}{\partial n} + \frac{\partial G_{|_{T_2}}}{\partial n}
\begin{theorem} 
Let $u$ be the solution of \eqref{2pBVP}. Then,  
\begin{equation}\label{2DGF}
\int_\Gamma u\, ds = \frac{1}{\sqrt{2}}  \int_\Omega  f G\, dx dy =  \frac{1}{\sqrt{2}}  \int_\Omega \nabla u \cdot \nabla G \, dx dy.
\end{equation}
\end{theorem}
\begin{remark} The simple calculations leading to formula \eqref{2DGF} have  the following  consequences:
\begin{enumerate}
\item [i)] The function $G=G(x,y)$ is the unique solution of 
\begin{equation}\label{2Ddelta}
  -\Delta u = \sqrt{2} \, \delta_\Gamma \ \text{on} \  \Omega =(0, 1)\times (0, 1), 
\end{equation}
where
\[
 \delta_\Gamma (\varphi) = \int_\Gamma \varphi\, ds, \ \text{for all} \ \varphi \in C_0^\infty (\Omega).
\]
\item [ii)] The function $G=G(x,y)$ is the Riesz representation of \\ the functional $ \sqrt{2} \, \delta_\Gamma :H^1_0(\Omega)\to \R$. 

\item [iii)] Part ii) can  be viewed  as an extension of  the 1D problem on $(0,1)$, where,  for each $s \in (0, 1)$, the function 
\[
x \to s(1-s) \Phi_s(x) = G(s,x)
\]
is the Riesz representation of  $\delta_s :H^1_0(\Omega)\to \R$, 
\[
\delta_s(\varphi) = \varphi(s), \ \text{for any} \ \varphi \in C_0^\infty (\Omega).
\]
\item [iv)]  The function  $G=G(x,y)$  can be viewed also as the finite element  discretization  of \eqref{2pBVP}, using $C^0-P^2$ discretization on the space $V_1$. 

\item [v)] In light of the FD versus FE connections presented here, any 2D quadrature approximation of the 2D integral 
$\displaystyle \frac{1}{\sqrt{2}}  \int_\Omega  f G$ can be viewed as a ``finite difference'' approximation of  $\displaystyle \int_\Gamma u\ ds$, where 
$u$ is  the solution of \eqref{2pBVP}.
\end{enumerate}
\end{remark}

%%%%%%%%%%%%%%
\subsection{ A bubble  function for the 2D discretization on a special domain}
%%%%%%%%%%%%%%
Consider that $\Omega$ is the domain defined by one equilateral triangle $T$ with vertices $z_1, z_2, z_3$. We define the bubble function
\[
B= \lambda_1 \lambda_2 \lambda_3,
\]
where $ \lambda_1, \lambda_2$ and  $\lambda_3$ are the  linear functions on $T$ with the property that 
$\lambda_i(z_j)=\delta_{ij},\  i, j=1,2,3$. We denote the area of the triangle $T$ by $|T|$. 

It is easy to check that $B\in H_0^1(\Omega)$ and 
\begin{equation}\label{B}
-\Delta B= \frac{1}{\sqrt{3}\, |T|}.
\end{equation} 
Then, for any function $\varphi \in H_0^1(\Omega)$, we have
\begin{equation}\label{B2}
 \frac{1}{\sqrt{3}\, |T|} \int_\Omega \varphi = \int_\Omega  (-\Delta B)  \varphi = \int_\Omega \nabla B\cdot \nabla \varphi .
\end{equation} 
In particular, if $u\in H_0^1(\Omega)$ is the solution of \eqref{2DBVP} for a given $f$, we have 
\begin{equation}\label{B3}
 \frac{1}{\sqrt{3}\, |T|} \int_\Omega u  = \int_\Omega \nabla B\cdot \nabla u =  \int_\Omega f\, B = (f, B).
\end{equation} 

\begin{theorem} 
Let $\Omega$ be an equilateral triangle, and let $u$ be the solution of \eqref{2pBVP}, then 
\begin{equation}\label{2DBF}
\int_\Omega u = \sqrt{3}|\Omega|  \int_\Omega  f B.
\end{equation}
\end{theorem}

\begin{remark} As in the previous case, we point out the following consequences:
\begin{enumerate}
\item [i)] The function $\sqrt{3}\, |\Omega|\, B$ is the unique solution of 
\[
 \left\{
 \begin{array}{rccl}
     -\Delta u = & 1  & \text{in } \ \Omega,\\
      u =& 0 & \text{on} \  \partial \Omega.
\end{array}  
\right. 
\]
\item [ii)]  From \eqref{B2},  the function $\sqrt{3}\, B$ is the Riesz representation of the functional $ F :H^1_0(\Omega)\to \R$,
\[
F(\varphi)= \frac{1}{ |\Omega|}  \int_\Omega \varphi.
\]

\item [iii)]  The function   $\sqrt{3}\, |\Omega|\, B$  can be viewed also as the  finite element  discretization  of \eqref{2pBVP}, using $C^0-P^3$ discretization on the space  \\ $V_0:= span\{B\}$. 
 \item [iv)] Any 2D quadrature approximation of 
$\sqrt{3} \int_\Omega  f B\, dx dy$ can be viewed as a ``finite difference'' approximation of  $\frac{1}{|\Omega|} \int_\Omega u$, where $u$ is  the solution of \eqref{2pBVP}.

\end{enumerate}
\end{remark}
We note that the assumption that $\Omega$ be an equilateral triangle is needed in order to have that the Laplacian operator acting on the bubble function is  a constant function.
%%%%%%%%%%%%%%
\section{Conclusion}\label{sec:conclusion}
%%%%%%%%%%%%%%

We considered the finite difference and the  finite element methods for a model  boundary value problem. 
We  emphasized the connections between the two methods by rewriting the FD system  such that the matrix of the system coincides with the stiffness matrix of the FE discretization.  In this reformulation, the right hand side vector for the FD system can be also viewed 
as a dual vector obtained componentwise from the composite trapezoid rule approximation of the FE dual vector.  
Using the connection between the Green function  and the $C^0-P^1$  finite element basis  functions, we  found that the inverse of the stiffness matrix $S$ is exactly the matrix obtained by evaluating the  Green function at the interior grid nodes $\{(x_i, x_j)\}  \subset (0, 1)  \times (0, 1)$.  Consequently, we found simplified proofs for the standard FD error analysis, and  provided an energy estimate
 for the difference between the FE and FD solutions.  
 We presented  interesting possible  extensions to the 2D case based on the fact that 
  the components of the FE  dual vector associated with the given data $f$  acting on  conforming finite element basis functions, can be written as linear combinations of various degrees of freedom acting on the solution $u$.  The challenges of extending these ideas to the more general cases of discrete spaces and meshes in 2D or 3D, remain to be addressed  in our future work.

\vspace{0.4in}

%%%%%%%%%%%%%%
\paragraph*{\bf Acknowledgment.}
%%%%%%%%%%%%%%
 The   authors are grateful for an Oden Institute fellowship which allowed them to visit UT Austin in the Fall of 2019.

\bibliographystyle{plain}

\end{document}